\documentclass[a4paper,twoside]{article}
\usepackage[latin1]{inputenc} 
\usepackage[T1]{fontenc} 
\usepackage{RR}
\usepackage{hyperref}
 \usepackage{graphics}
 \usepackage{graphicx}
 \usepackage{epsfig}
\usepackage{amssymb}
\usepackage{amsthm}
\usepackage{amsmath}
\usepackage{url}
\usepackage{makeidx}
\usepackage{multicol}
\usepackage{psfrag}
\usepackage{color}

\RRNo{7803}
\RRdate{November 2011}
\RRversion{2}
 \RRdater{January 2014}
\RRauthor{
Ibtihel {\sc Ben Gharbia}\thanks{Eddimo, Ceremade, University of
Paris-Dauphine, 75775~Paris, France}%
 \thanks[INRIA]{INRIA-Paris-Rocquencourt, team-project Pomdapi, BP~105,
 F-78153~Le Chesnay, France  
 e-mails\,:  \texttt{Ibtihel.Ben\_Gharbia@inria.fr};  \texttt{Jerome.Jaffre@inria.fr}
 }%
  \and
J\'er\^ome {\sc Jaffr\'e}
\thanksref{INRIA}
}
\authorhead{I. Ben Gharbia \& J. Jaffr\'e}
\RRtitle{Conditions de compl\'ementarit\'e pour l'apparition et la disparition de la phase gazeuse }
\RRetitle{Gas phase appearance and disappearance as a problem with complementarity constraints}
\titlehead{Complementarity constraints for gas phase appearance and disappearance }
\RRnote{This work was partially supported by the GNR MoMaS (PACEN/CNRS, ANDRA, BRGM, CEA, EDF, IRSN)}
\RRresume{La migration d'hydrog\`ene produit par la corrosion des sites de stockages souterrains
des d\'echets nucl\'eaires avec dissolution de l'hydrog\`ene est formul\'ee comme un ensemble  d'\'equations aux d\'eriv\'ees partielles non-lin\'eaires
avec des conditions de compl\'ementarit\'e non-lin\'eaires. Cet article montre comment appliquer une strat\'egie moderne et efficace, la m\'ethode de Newton-min,  pour r\'esoudre ce  probl\`eme de g\'eosciences. En particulier, les exp\'eriences num\'eriques montrent que la m\'ethode de Newton-min se r\'ev\'ele efficace et converge quadratiquement pour ce probl\`eme. }
\RRabstract{
The modeling of migration of hydrogen produced by the corrosion of the nuclear waste packages in 
an underground storage including the dissolution of hydrogen involves a set of nonlinear partial differential
 equations with nonlinear complementarity constraints. This article shows  how to apply a modern 
 and efficient solution strategy, the Newton-min method, to this  geoscience problem and  investigates 
 its applicability and efficiency.  In particular, numerical experiments show that the Newton-min method 
 is quadratically convergent for this problem. }
\RRmotcle{Milieu poreux,
\'ecoulement diphasique, 
dissolution, 
stockage profond de d\'echets nucl\'eaires,
probl\`eme de compl\'ementarit\'e non-lin\'eaire,
fonction non-lisse,  
 Newton-min
}
\RRkeyword{Porous media,  
two-phase flow, 
dissolution, 
nuclear waste underground storage,
 nonlinear complementarity problem, 
non-smooth function, 
Newton-min}
\RRprojets{Pomdapi}
 \URRocq 
 \RCParis 

\newcommand{\RR}{{\mathbb R}}
\newcommand{\T}{^{\mskip-1mu\top\mskip-2mu}}

\let\leq\leqslant
\let\geq\geqslant
\if{

}\fi

\begin{document}
\makeRR   
\tableofcontents
\section{Introduction}
\label{Introduction}
The couplex-Gas benchmark~\cite{couplex} was proposed by Andra {(French National Inventory of Radioactive Materials and Waste)~\cite{andra}} and {the research group MoMaS (Mathematical Modeling and Numerical Simulation for Nuclear Waste Management Problems)~\cite{momas}} in order to improve
 the simulation of the migration of hydrogen produced by the corrosion of nuclear waste packages in 
 an underground storage. This is a system of two-phase (liquid-gas) flow with two components (hydrogen-water). 
 The benchmark generated some interest and engineers encountered difficulties in handling the appearance 
 and disappearance of the phases. The resulting formulation \cite{JS09} is a set of partial differential equations 
 with nonlinear complementarity constraints. Even though they appear in several problems of flow and transport 
 in porous media like the black oil model presented in \cite{CJ1986} or transport problems with dissolution-precipitation
 \cite{Krautle08,BKKK2011,MMK1}, complementarity problems are not usually identified as such in hydrogeology and, 
 to circumvent the solution of complementarity conditions,  problems are often solved by reformulating the problem
 as in~\cite{BJS09,Abadpour-Panfilov-2010,AAPP2010}. However the solution of complementarity problems is an active 
 field in optimization \cite{BGG10,facchinei-pang-2003,harker-pang-1990} and we draw from the know-how of this scientific 
 community. A similar path is followed bin papers like~\cite{MMK2,hw-2010,lhhw-2011}. 
The application of a semi-smooth Newton method~\cite{hintermueller-ito-kunisch-2003,kanzow-2004}, sometimes called the Newton-min algorithm, to solve nonlinear
complementarity problem is described. We will demonstrate through a test case, the ability of our model and our
solver to efficiently cope with appearance or/and disappearance of one phase.

In the section \ref{probform}, we introduce the formulation of the problem and in the section \ref{nummeth} 
we describe the numerical method.  In the section~\ref{numsim}, we present and discuss a numerical experiment.

\section{Problem formulation}
\label{probform}
This section gives a precise formulation of the mathematical model for the application that was outlined in the
introduction. We consider a problem where the gas phase can disappear while the liquid phase is always present. 
\subsection{Fluid phases}
Let $\ell$ and $g$ be the respective indices for the liquid phase and the gas phase. 
Darcy's law reads
\begin{equation}
\label{darcy}
{\bf q}_i = -K(x) k_i(s_i) (\nabla p_i - \rho_i g \nabla z), \quad i=\ell, g,
\end{equation}
where $K$ is the absolute permeability. For each phase $i = \ell,\, g$, $s_i$ is the saturation 
and $k_i = \dfrac{{k_r}_i(s_i)}{\mu_i}$ is the mobility with $k_{ri}$ the relative permeability 
and $\mu_i$ the viscosity (assumed to be constant). 
The mobility  $k_i$ is an increasing function of $s_i$ such that $k_i(0) = 0, \; i = \ell,g$.
Assuming that the phases occupy the whole pore space, the phase saturations satisfy 
\begin{center}
$0 \leq s_i \leq 1, \quad s_\ell+s_g=1$.
\end{center}
The phase pressures are related through the capillary pressure law 
\begin{center}
$p_c(s_\ell) = p_g -p_\ell  \geq 0$, 
\end{center}
assuming that the gas phase is the non-wetting phase. 
The capillary pressure is a decreasing function of the saturation $s_\ell$. 

In the following, we will choose $s_\ell$ and $p_\ell$ as the main variables since we assume that 
the liquid phase cannot disappear for the problem under consideration. 

\subsection{Fluid components}
\medskip
We consider two components, water and hydrogen, identified by the indices $j=w$,~$h$. 
The mass density of the phase is
$$
\rho_i = \rho_w^i + \rho_h^i, \quad i=\ell,g.
$$
From $M^w$ and $M^h$, the water and hydrogen molar masses, we define the molar concentration of phase $i$:
\begin{equation}
c_i = c_w^i + c_h^i, \quad c_j^i = \dfrac{s_i\rho_j^i}{M^j},\quad j = w,h, \quad  i=\ell,g.
\end{equation}
The molar fractions are
\begin{equation}
\chi_h^i = \frac{c_h^i}{c_i},\quad \chi_w^i = \frac{c_w^i}{c_i}, \quad  i=\ell,g.
\end{equation}
Obviously,  
\begin{equation}
\chi_w^i + \chi_h^i = 1,\quad  i=\ell,g.
\end{equation}
We assume that the liquid phase may contain both components, while the gas phase contains only hydrogen, 
that is the water does not vaporize. In this situation we have
\begin{center}
$\rho_w^g=0, \quad \rho_g = \rho_h^g, \quad \chi_h^g = \frac{c_h^g}{c_g} =1 ,\quad \chi_w^g =0.$
\end{center}
{For the liquid phase, we assume that the water is the solvent and the hydrogen is the solute and that the quantity of hydrogen dissolved in the liquid is small, that is $ c_h^l \ll  c_w^l$. So we have
$$
 \chi_h^\ell \approx \dfrac{c_h^\ell}{c_w^\ell} = \dfrac {M^w}{M^h\rho_w^\ell}\rho_h^\ell. 
 $$
 }
 A third main unknown will be $\chi^\ell_h$, in addition to $s_\ell$ and $p_\ell$.
\subsection{Conservation of mass}
\medskip
We introduce the molecular diffusion flux for the diffusion of hydrogen in the liquid phase
\begin{equation}
j_h^\ell = -\phi M^h s_\ell c_\ell D_h^\ell\nabla\chi_h^\ell
\end{equation}
where $D_h^\ell$  is a molecular diffusion coefficient.

Conservation of  mass applied to each component, water and hydrogen, gives 
\begin{equation} \begin{array}{l}
\label{cons-equa}
 \dfrac{\partial}{\partial t} (\phi{\rho_{w}^\ell} s_\ell) + 
 \mbox{div} ({\rho_{w}^\ell}{\bf q}_\ell -  j_h^\ell ) = Q_w,\\[0.3cm]
  \dfrac{\partial}{\partial t} (\phi s_\ell \rho^\ell_h +\phi s_g\rho^g_h) +
   \mbox{div}(\rho^\ell_h {\bf q}_\ell + \rho^g_h {\bf q}_g + j_h^\ell) = Q_h.
\end{array} \end{equation}
We assume also that the gas is slightly compressible, that is $\textstyle{\rho_g = C_g p_g}$ with $C_g$ 
the compressibility constant, and that the liquid phase is incompressible, that is $\rho^\ell_w$ is constant.  
\subsection{Nonlinear complementarity constraints}
\medskip
Next, we apply Henry's law which {says} that, at a constant temperature, the amount of a given gas 
that dissolves in a given type and volume of liquid is directly proportional to the partial pressure 
of that gas in equilibrium with that liquid.

In the presence of the gas phase, Henry's law reads $H p_g = \rho_h^\ell,$ where $H= \mbox{H}(T) M^h$ 
with H(T) is the Henry law constant, depending only on the temperature.

There are two possible cases\,: the gas phase exists: $1-s_\ell > 0$, Henry's law applies 
and $H(p_\ell + p_c(s_\ell))  - \rho_h^\ell = 0$, or the gas phase does not exist, {$s_\ell = 1$
and $H(p_\ell + p_c(1))  - \rho_h^\ell \geq 0$} which says that for a given pressure $p_\ell$ the concentration $\rho_h^\ell$ 
is too small for the hydrogen component to be partly gaseous, or conversely for a given concentration 
$\rho_h^\ell$ the pressure $p_\ell$ is too large for the hydrogen component to be partly gaseous.

These cases can be written as complementary constraints
\begin{equation} \begin{array}{l}
\label{comp}
(1-{ s_\ell})\bigl(H({p_\ell} + p_c({ s_\ell})) - \rho_h^\ell\bigr) =0,~~
1-{ s_\ell} \geq 0,~~
H({p_\ell} + p_c({ s_\ell})) - \rho_h^\ell \geq 0.
\end{array} \end{equation}
Finally we end up with a system of nonlinear partial differential equations (conservation
equations \eqref{cons-equa} and Darcy laws \eqref{darcy}) 
with the nonlinear complementarity constraints \eqref{comp} describing the transfer
of hydrogen between the two phases, the unknowns being $s_\ell$, $p_\ell$, and $\chi_h^\ell$. 
This formulation has the advantage of being valid whether the gas phase exists or not \cite{JS09}.

\section{ Discretization and solution method}
\label{nummeth}
We {use} a first order Euler implicit scheme for time discretization and cell-centered finite volumes 
for space discretization.
We {denote} by $N$, the number of degrees of freedom for $s_\ell$, $p_\ell$ and $\chi^\ell_h$ 
which is equal to the number of cells. We introduce
\begin{itemize}
\item ${x} \in \RR^{3N}$, the vector of unknowns for ${s_\ell}$, ${p_\ell}$, ${\chi_h^\ell}$,
\item ${{\cal H}}: \RR^{3N} \rightarrow \RR^{2N}$, the discretized conservation equations,
\item ${{\cal F}}: \RR^{3N} \rightarrow \RR^{N}$, the discretized function $1-{s_\ell}$,
\item ${{\cal G}}: \RR^{3N} \rightarrow \RR^{N}$, the discretized function 
{$H({p_\ell} + p_c({s_\ell})) -  \dfrac{M^h\rho_w^\ell}{M^w} {\chi_h^\ell}$}.
\end{itemize}
Then at each time step the problem can be written in compact form
\begin{equation}
\begin{array}{l}
\label{comp-problem}
{{\cal H}(x) =0,}\\[0.2cm]
{{\cal F}(x) \T {\cal G}(x) =0,\quad {\cal F}(x) \geq 0, \quad {\cal G}(x) \geq 0,}
\end{array}
\end{equation}
where the inequalities have to be understood component-wise.

\subsection{A non-smooth system using the Minimum function}
\label{solution-method}
\medskip
It is well known that complementarity conditions, consisting of equations and inequalities, 
can be expressed equivalently by an equation via a complementarity function~\cite{facchinei-pang-2003}(C-function).
Let 
{
$$\begin{array}{ccccc}
\varphi & : & \RR^N \times\RR^N  & \to &\RR^N \\
& & (a,b) & \mapsto & \min(a,b)\\
\end{array}
$$
}
be the minimum function, in which the $\min$ operator acts component-wise. 
This is a C-function, in the sense that it satisfies
\begin{equation}
\label{complementarity1}
{\varphi(a,b)} =0\qquad \Longleftrightarrow\qquad  a \geq 0, \quad b \geq 0,\quad a\T b =0.
\end{equation}
Other typical scalar C-functions~\cite{facchinei-pang-2003} are
\begin{itemize}
\item the Fisher-Burmeister function\,: $\varphi(a,b)=\sqrt{a^2+b^2}-a-b$, 
\item $\varphi(a,b)=-ab+\min^2(0,a)+\min^2(0,b).$
\end{itemize}
Using this minimum function, we can write the complementarity problem \eqref{comp-problem} as
\begin{equation}
 \begin{array}{l}
 \label{prob}
{{\cal H}(x) =0,}\\[0.2cm]
\varphi({\cal F}(x), {\cal G}(x)) = 0.
\end{array}
\end{equation}
Hence, the resulting system of mass conservation (differential) equations and equilibrium conditions 
is fully free of inequalities (pure set of equations).
The only drawback of the introduction of a complementarity problem is that the problem is no longer $C^1$,
since $\varphi \notin C^1(\RR^{2N},\RR^N)$, while the typical assumption for having the local quadratic 
convergence of Newton's algorithm requires to have a ``$C^1$ function with a Lipchitz-continuous derivative".
However, it is well known, especially in the community of optimization, 
that the assumptions can be weakened in several ways, for example by only assuming strong semi-smoothness.
In the next section we give the definition of semi-smoothness from~\cite{clarke-1990,facchinei-pang-2003}.

\subsection{Semi-smoothness}
\label{semi-smoothness}

\medskip\noindent
Let $\psi : \RR^N \to \RR^N$ be a locally {lipschitz-continuous} function. Then, by Rademacher's
theorem\cite{facchinei-pang-2003}, there is a dense subset $D \subset \RR^N$ on which $f$ is differentiable.
The $B$-subdifferential of $\psi$ at a point~$x\in\RR^N$ is the set
$$
\partial_B \psi(x):=\{ J \in \RR^{N\times N} \quad |  \quad J = \lim_{k\to\infty}
\psi^{\prime}(x_k),~(x_k) \subset D,~ x_k \to x\},
$$
where $\psi^{\prime}$ is the derivative of $\psi$.
The {\em
generalized Jacobian\/} of ${\psi}$ at~$x$~\cite{clarke-1990} is the set
$$
\partial \psi(x)=\mbox{co} \,\partial_B \psi(x),
$$
where $\mbox{co}\, S$ denotes the convex hull of a
set $S$.
Now, the function
$\psi$ is said to be semi-smooth at $x$ if {$\psi$ is directionally differentiable} at $x$ and
\begin{center}
$Jd - \psi^{\prime}(x;d) = o(||d||)$,
\end{center}
for any $d\to0$ and for any $J \in \partial\psi(x+d)$, where $\psi^{\prime}(x;d)$ denotes the directional
derivative of $\psi$ at $x$ in the direction of $d$. Analogously, $\psi$ is called {\em strongly
semi-smooth\/} at~$x$, if 
$$
Jd - \psi^{\prime}(x;d) = o(||d||^2).
$$
$\psi$ is called (strongly) semi-smooth if $\psi$ is (strongly) semi-smooth at any point~{$x~\in~\RR^N$}.
\medskip

 It is well known that the minimum function and  the Fisher-Burmeister function are strongly semi-smooth. One
 can then solve system \eqref{prob} using the semi-smooth Newton's method,  called 
 the Newton-min method \cite{BGG10,bengharbia-gilbert-2012} when the min function is used. The Newton-min method can also be regarded as 
 an active set strategy \cite{hintermueller-ito-kunisch-2003}.

\subsection{The Newton-min algorithm}
\label{algorithm}
\medskip
We now give an exact statement of the Newton-min algorithm for solving the nonlinear 
system of equation (\ref{prob}).

Below $\partial {\cal \varphi} (x)$ denotes the generalized Jacobian of ${\cal \varphi}$ at a point~$x$.
Let Res be the residual of~${\cal \psi}(x)$~where~${ {\cal \psi}(x) : =
\begin{pmatrix}
{{\cal H}(x)}\\
{\cal \varphi}(x)\\
\end{pmatrix}}
$
and $\varepsilon$ be a stopping criterion for Res.
\medskip
\begin{quote}
\rule{\linewidth}{0.5pt}
\\[2ex]
\medskip
\noindent
Let {$x^1\in\RR^N$}.
For $k=2,3,\ldots$, do the following.
\begin{list}{}{\topsep=0.5ex\parsep=0.0ex\itemsep=0.0ex
 \settowidth{\labelwidth}{9)}
 \labelsep=1ex
 \leftmargin=\labelwidth\addtolength{\leftmargin}{\labelsep}\addtolength{\leftmargin}{1ex}
}
\item[1)]
\label{step1}
If $\mbox{Res} \leq \varepsilon$, stop.
\medskip
\item[2)]
\label{step2}
Define the complementary index sets $A^k$ and
$I^k$ by
$$
A^k:=\{i:{\cal G}_i(x^{k})<{\cal F}_i(x^{k})\}, \quad
I^k:=\{i:{\cal G}_i(x^{k}) \geq {\cal F}_i(x^{k})\}.
$$
\item[3)]
\label{step3}
Select an element ${\cal J}_x^k \in \partial {\cal \varphi} (x^k)$ such that 
 its $i$th line is equal to ${\cal
F}_i'(x^k)$ [resp.\ ${\cal G}_i'(x^k)$] if ${\cal F}_i(x^k)\leq {\cal G}_i(x^k)$ [resp.\ ${\cal F}_i(x^k)> {\cal
G}_i(x^k)$].
\medskip
\item[4)]
Let $x^{k+1}$ be a solution to
\[ \begin{array}{l}
{{\cal H}(x^k) + {\cal H}^{'}(x^k)(x^{k+1}-x^k) =0,}\\[0.1cm]
{{\cal \varphi}(x^k) + { {\cal J}_x^k} (x^{k+1}-x^k)  =0,}  \qquad {{\cal J}_x^k \in \partial {\cal \varphi} (x^k).}   \\[0.1cm]
  \end{array}\]
\end{list}
\vspace{2ex}
\rule{\linewidth}{0.5pt}
\end{quote}
\bigskip
\noindent

Note that, as in a smooth Newton method, only one linear system has to be solved at each Newton iteration.

Furthermore the Newton-min method satisfies also a quadratic convergence property. Indeed, a theorem\cite{facchinei-pang-2003} says that if  $x^{*}$ is
a solution to the system $\psi(x) = 0$, such that $J$ is nonsingular for all $J \in \partial \psi(x^*)$  (as defined in section~\ref{semi-smoothness}), then for any initial value sufficiently close to $x^{*}$, the Newton-min method generates a sequence
that converges quadratically to $x^*$. 

We have not yet proved the hypothesis of non-singularity of $J$ for our system but we observed the quadratic convergence in our numerical experiments.

\section{Numerical experiment}
\label{numsim}
\subsection{A problem inspired from the Couplex Gas benchmark}
\medskip
We consider a one-dimensional core with length $\mbox{L}=200\,\mbox{m}$, initially saturated with liquid $(s_\ell =1)$
and containing no hydrogen {$(\chi_h^\ell =0)$}. 
Hydrogen is injected at a given rate on the left. After a while the hydrogen injection is stopped. 
The problem is then to simulate the migration of hydrogen and to illustrate the gas appearance and disappearance phenomena.

We calculate {spatial} evolutions of the liquid pressure, the total hydrogen molar density and the the gas saturation along the line. Computations are performed from the initial time up to the stationary state.

The core is supposed to be homogenous porous medium. The capillary pressure function $p_c$ 
and the relative permeability functions, $k_{rl}$ and $k_{rg}$, are given by the Van Genuchten-Mualem model~\cite{VG1980}:
$$
\begin{array}{ll}
 p_c      =  P_r \left(S^{-{1}/{m}}_{le}-1\right)^{1/n} ,&\\
 k_{rl}    = \sqrt{S_{le}}\left(1-\left(1-S^{{1}/{m}}_{le}\right)^m\right)^2, & k_{rg} = \sqrt{1-S_{le}}\left(1-S^{{1}/{m}}_{le}\right)^{2m},
\end{array}
$$
with {$S_{le} = \dfrac{S_l-S_{lr}}{1-S_{lr}-S_{gr}}$} and $m = 1 - \dfrac{1}{n}$, and where parameters $P_r$, $n$, $S_{lr}$ and $S_{gr}$ depend on the porous medium. {The parameters} describing the porous medium and {the fluid} characteristics are given in Table \ref{tab1}. Fluid temperature is fixed to $T = 303$ K.
\begin{table}[hbtp]
\begin{center}
\begin{tabular}{|c|c|c|c|} 
\hline
\multicolumn{2}{|c|}{Porous medium parameters} & \multicolumn{2}{|c|}{Fluid characteristics parameters} \\
\hline
 Parameter \:&  Value  & Parameter\;&  Value  \\
  \hline
      $K$             & 5\,10$^{-20}$  \; m$^2$  & $T$                                     & 303 \quad K\\
      $\phi$        &  0.15 \quad (-)                   & $D^h_\ell$                        & 3 10$^{-9}$ \quad m$^2$/s \\
      $P_r$         &  2 10$^6$ \quad Pa        & $\mu_\ell$                         & {1 10$^{-9}$ \quad Pa.s}\\
      $n$             & 1.49 \quad (-)                   & $\mu_g$                        & {9 10$^{-9} $\quad Pa.s} \\ 
      $S_{lr}$     &  0.4 \quad(-)                      & $H(T = 303 \mbox{K})$  & 7.65 10$^{-6}$\, mol/Pa/m$^3$\\ 
      $S_{gr}$    &  0   \quad    (-)                   & $M_w$                              & 10$^{-2} $\quad kg/mol \\ 
                          &                                             & $M_h$                               & 2 10$^{-3} $\quad kg/mol \\
                          &                                             &  {$\rho^\ell_w$}                     & 10$^{3} $  \quad kg/m$^3$ \\ 
\hline
\end{tabular}
\end{center}
\caption{Values of porous medium fluid characteristics.  }
\label{tab1}
\end{table}  

Initial conditions are $S_\ell \,(t=0) = 1$, {$\chi_h^{\ell}\, (t=0) = 0$} and~$p_\ell\,(t=0) =10^6$~Pa.
For boundary conditions on the left, the hydrogen flow rate is given, 
$\rho^\ell_h {\bf q}_\ell + \rho^g_h {\bf q}_g + j_h^\ell = 5.57\,10^{-6}$ kg/m$^{2}$/year. 
From this condition, one can deduce the saturation. 
Still on the left, we impose a zero water flow rate ${\rho_{w}^\ell}{\bf q}_\ell -  j_h^\ell = 0$. 
On the right, the liquid pressure is given, $p_\ell = 10^6$ Pa, and the liquid saturation {is set to $s_\ell =1$}. 

\subsection{Results and comments}
\medskip
For the numerical simulation below we divided the space interval into 200 intervals of equal length and we used
a constant time step of 5000 years.
During the simulation, we can identify four important periods, three periods during injection 
and one period after injection.

\begin{figure}
\begin{center}
\includegraphics[width=9cm,height=7cm]{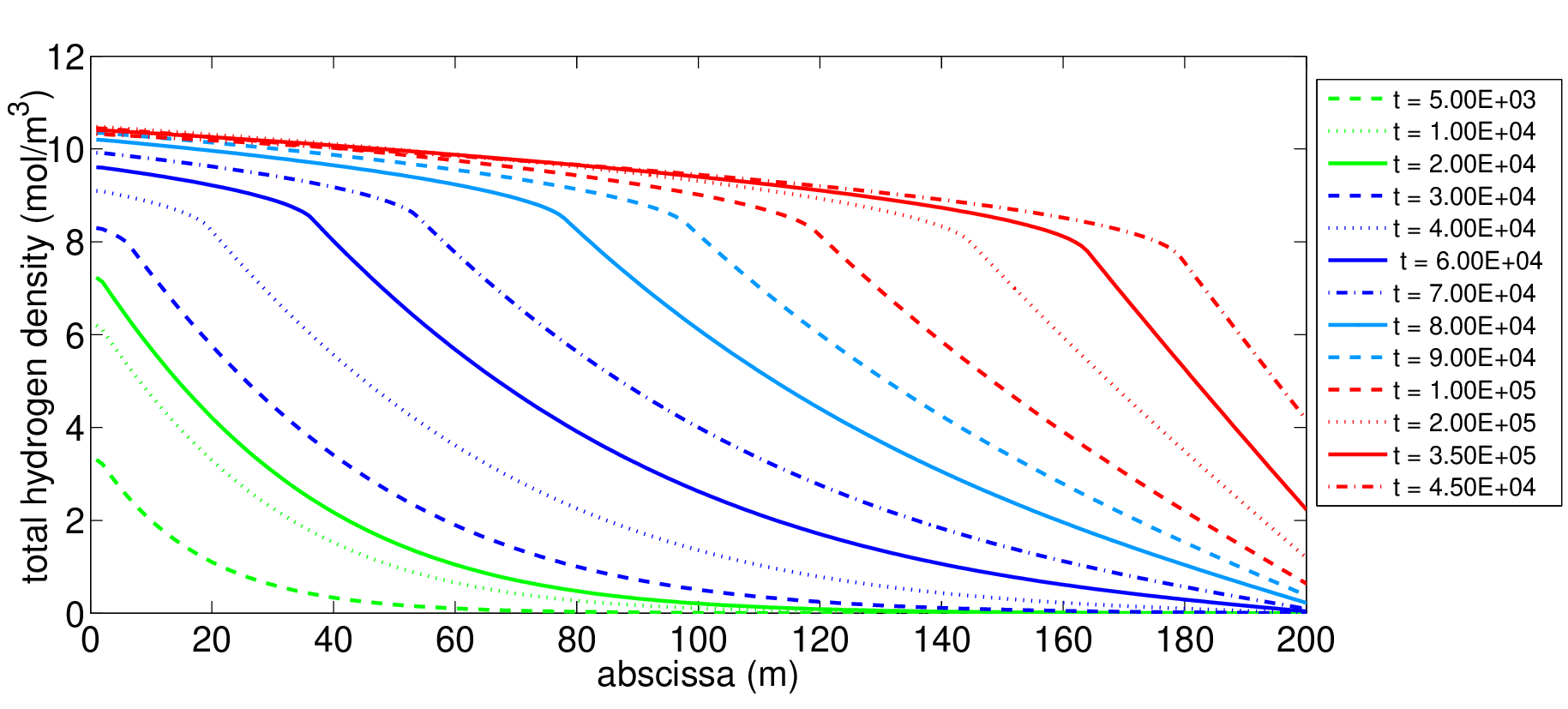}
\caption{\label{during-injection-rho}
Spatial evolution of  hydrogen density at several times t (in years) during hydrogen injection.}
\end{center}
\end{figure}

\begin{figure}
\begin{center}
\includegraphics[width=10cm,height=7cm]{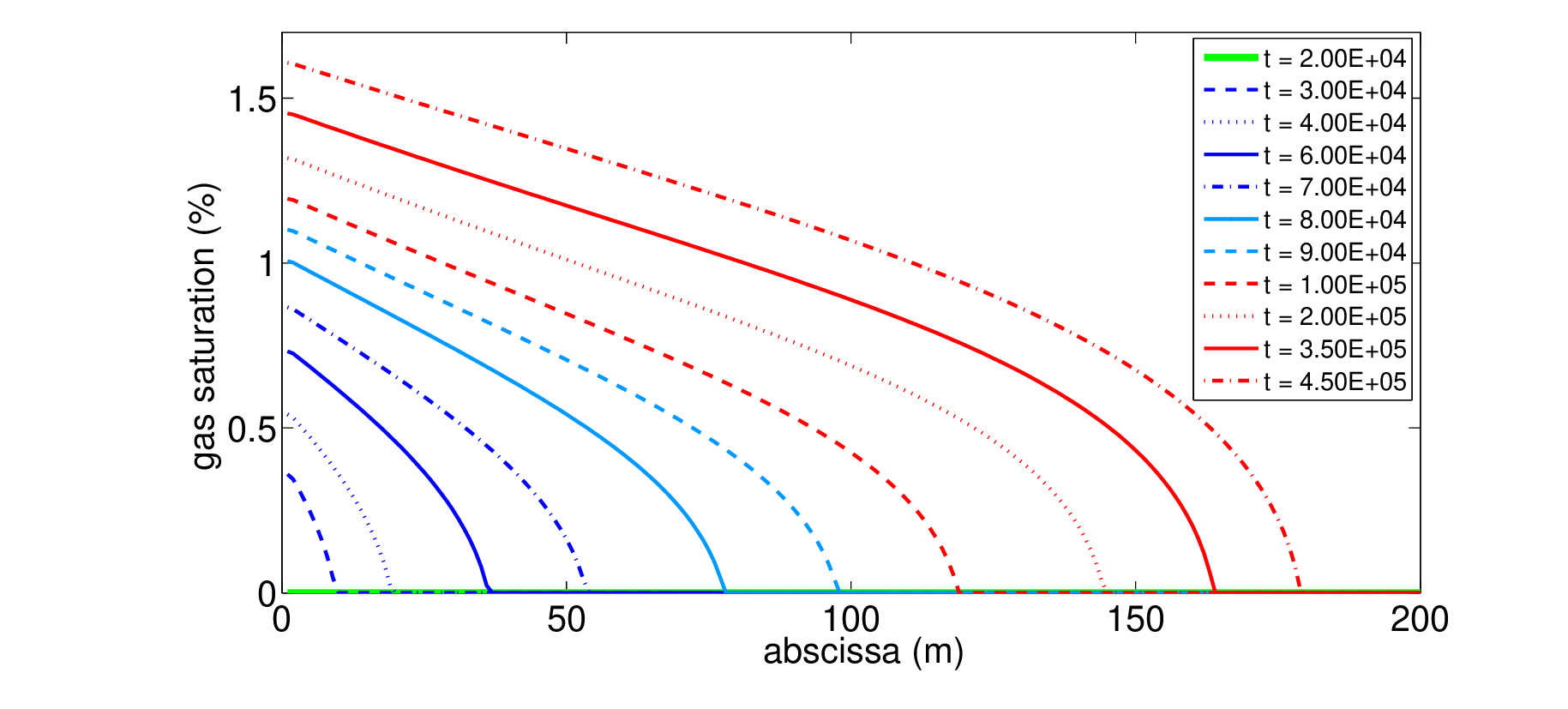}
\caption{\label{during-injection-s}
Spatial evolution of gas saturation at several times t (in years) during hydrogen injection.}
\end{center}
\end{figure}

\begin{figure}
\begin{center}
\includegraphics[width=9cm,height=7cm]{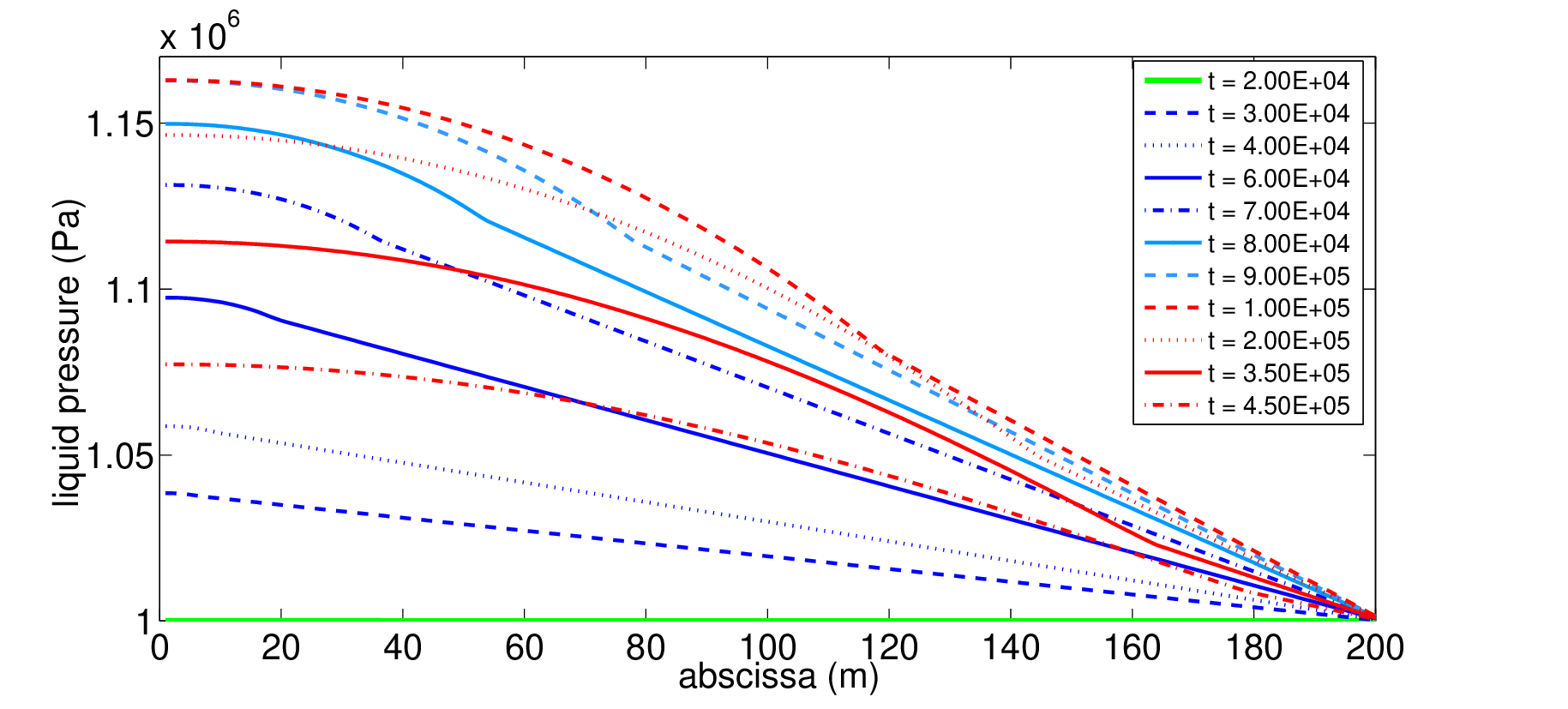}
\caption{\label{during-injection-p}
Spatial evolution of liquid pressure at several times t (in years) during hydrogen injection.}
\end{center}
\end{figure}

\begin{figure}
\begin{center}
\includegraphics[width=9cm,height=7cm]{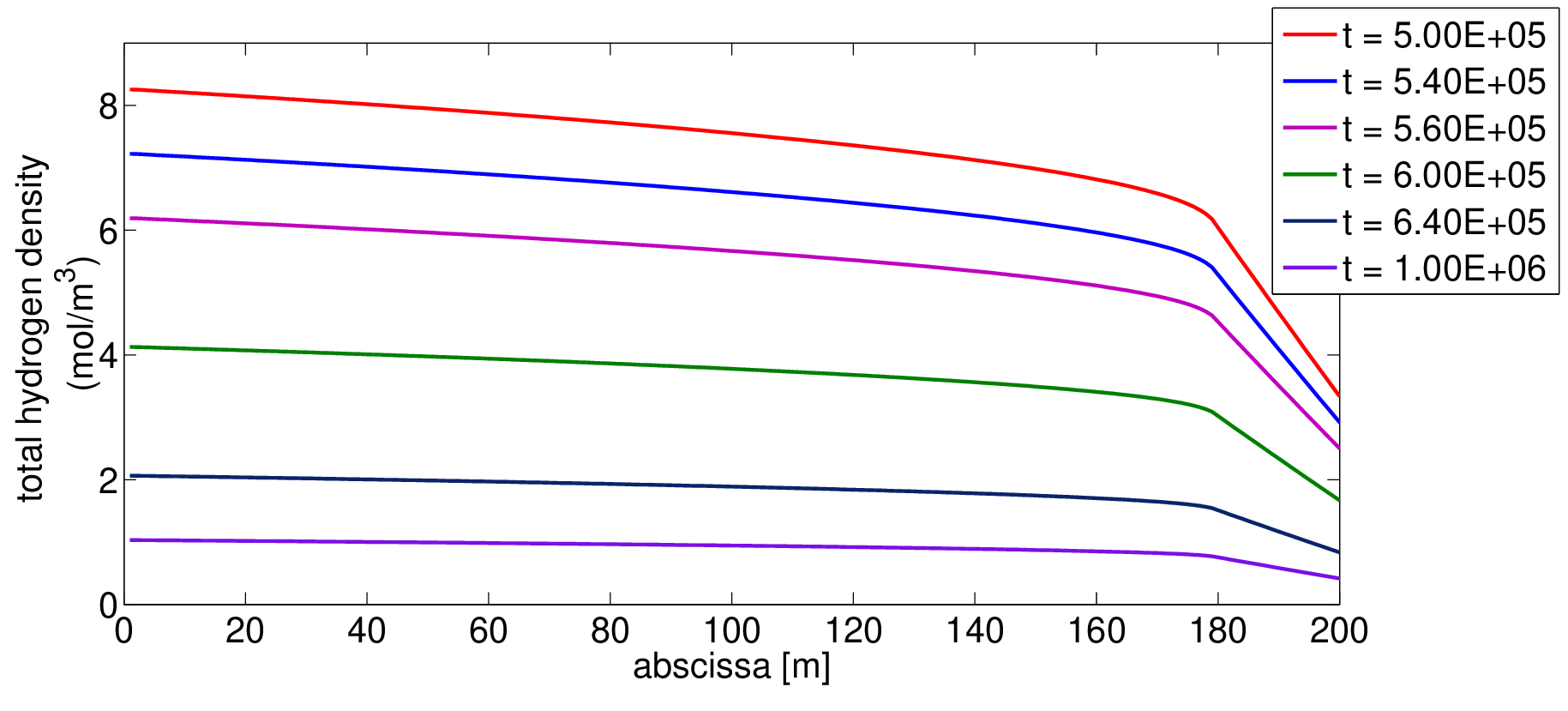}
\caption{\label{after-injection-rho}
Spatial evolution of  hydrogen density at several times t (in years) after hydrogen injection is stopped.}
\end{center}
\end{figure}

\begin{figure}
\begin{center}
\includegraphics[width=9cm,height=7cm]{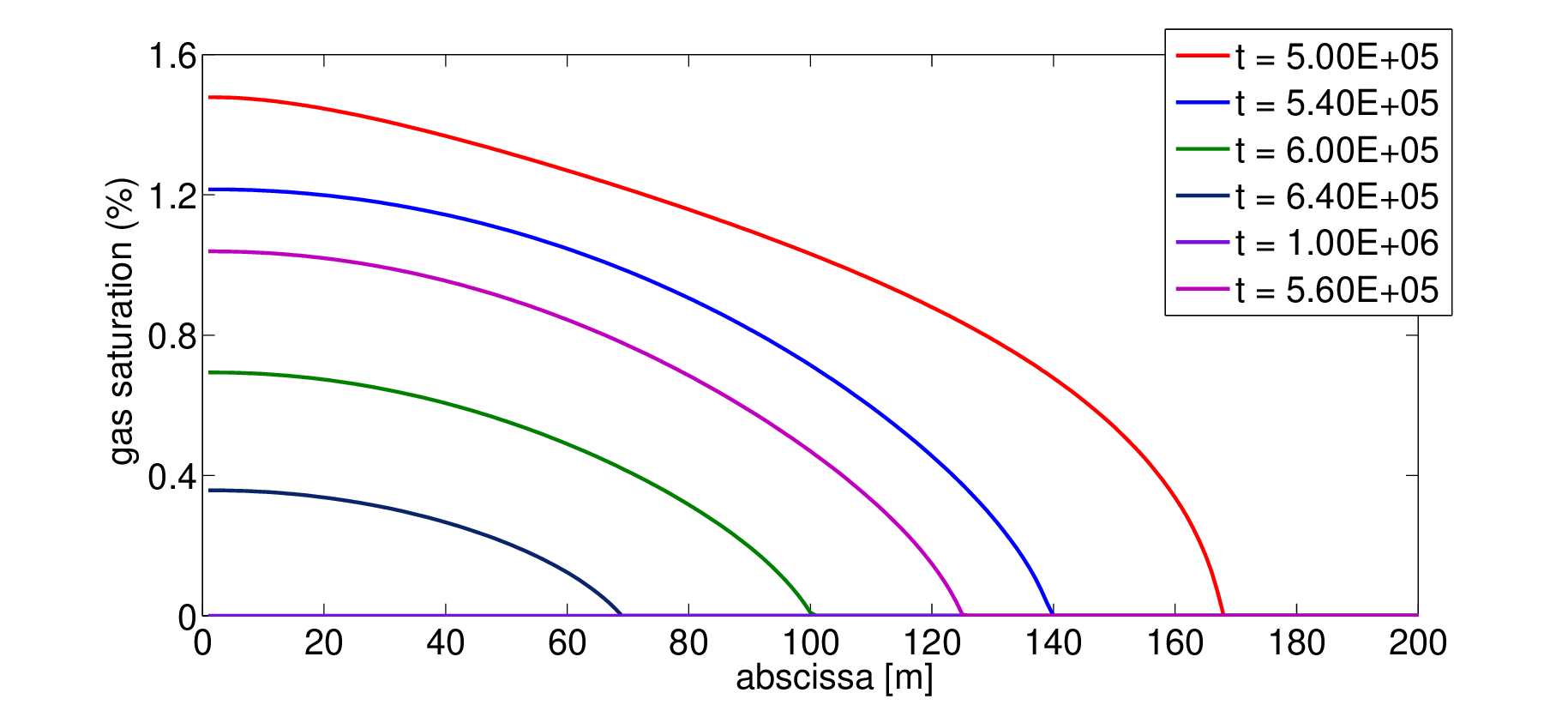}
\caption{\label{after-injection-s}
Spatial evolution of gas saturation at several times t (in years) after hydrogen injection is stopped.}
\end{center}
\end{figure}

\begin{figure}
\begin{center}
\includegraphics[width=9cm,height=7cm]{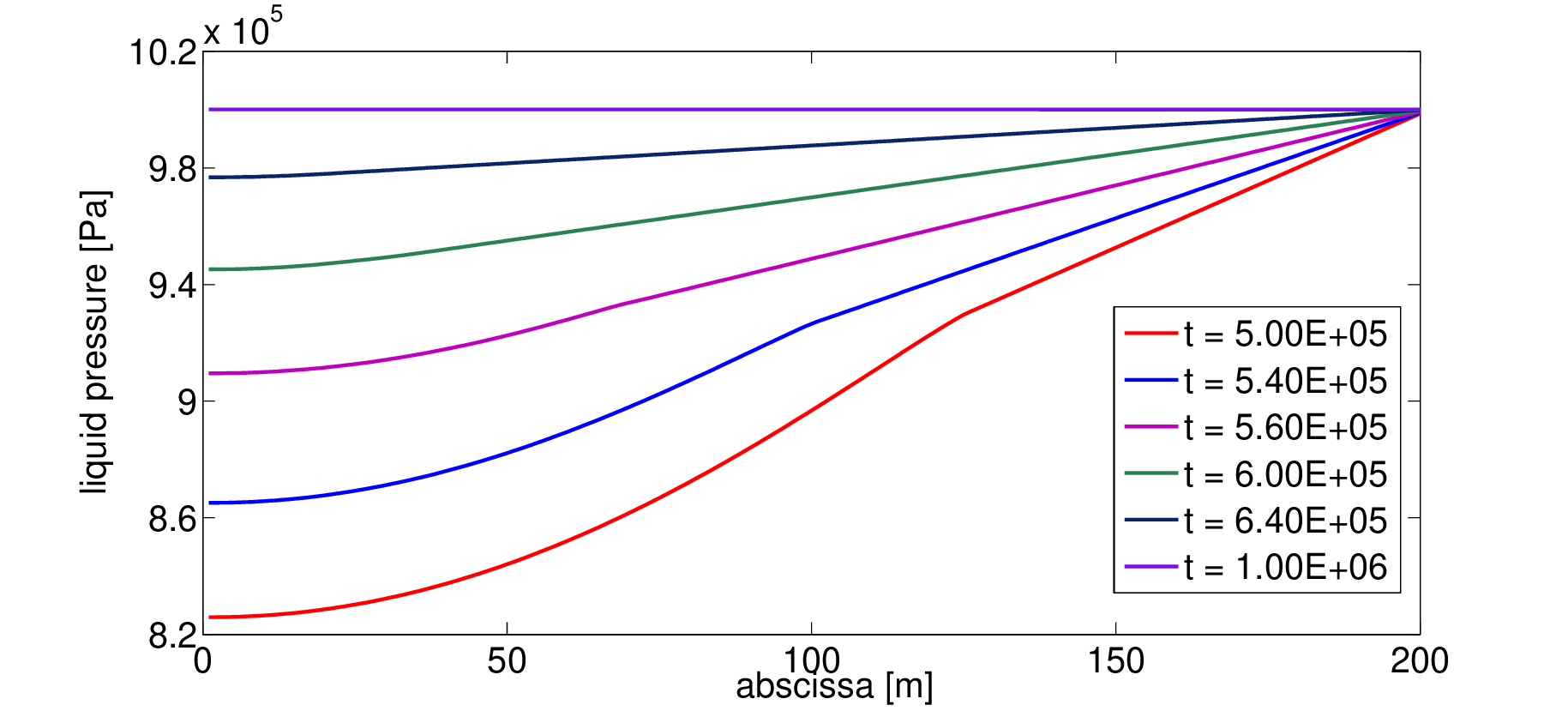}
\caption{\label{after-injection-p}
Spatial evolution of liquid pressure at several times t (in years) after hydrogen injection is stopped.}
\end{center}
\end{figure}

\bigskip
\noindent
{\bf  During injection} (figures~\ref{during-injection-rho}, \ref{during-injection-s} and \ref{during-injection-p}): ${\bf0<t<5.10^5}$ years
\begin{itemize}
\item Period $1$ (${\bf0 < t <2 \;10^4}$ years): only the hydrogen density increases 
(Figure \ref{during-injection-rho}, green curves),
while the liquid pressure and the gas saturation stay constant (Figures \ref{during-injection-s} and~\ref{during-injection-p}, green curve); the whole domain is saturated with water $(s_g=0)$.

\item Period $2$ (${\bf 2 \;10^4 \leq t \leq 1.5 \;10^5} $ years): at $t = 2 \;10^4$, the gas phase appears ($s_g > 0$). 
During this period, the liquid pressure increases (Figures~\ref{during-injection-p}, blue curves) and pressure gradients are non zero which corresponds to a displacement
of both phases. 
The total hydrogen density and the gas saturation increase (Figures~\ref{during-injection-rho} and \ref{during-injection-s}, blue curves) and the unsatured 
area grows.

\item Period $3$ ($\bf 1.5 \;10^5 < t < 5 \;10^5 $ years):  while the total hydrogen density 
and the gas saturation continue to increase (Figures \ref{during-injection-rho} and~\ref{during-injection-s}, red curves); the liquid pressure and the pressure gradient
 decrease since there is no water injection (Figure~\ref{during-injection-p}, red curves).\\[0.05cm]
\end{itemize}
\noindent
{\bf After injection} (Figures~\ref{after-injection-rho},~\ref{after-injection-s} and~~\ref{after-injection-p}): 
\begin{itemize}
\item Period $4$ (${\bf t > 5 \;10^5}$ years): cell by cell, starting from the right, the gas saturation decreases and after a while, the gas phase disappears (Figure~\ref{after-injection-s}).  At the end of the simulation the system reaches a stationary state (Figure~\ref{after-injection-rho}) and the liquid pressure gradient goes to zero (Figure~\ref{after-injection-p}).
\end{itemize}

\newpage
\subsection{Quadratic convergence}
\begin{center}
\begin{figure}[h]
\includegraphics[width=15cm,height=9cm]{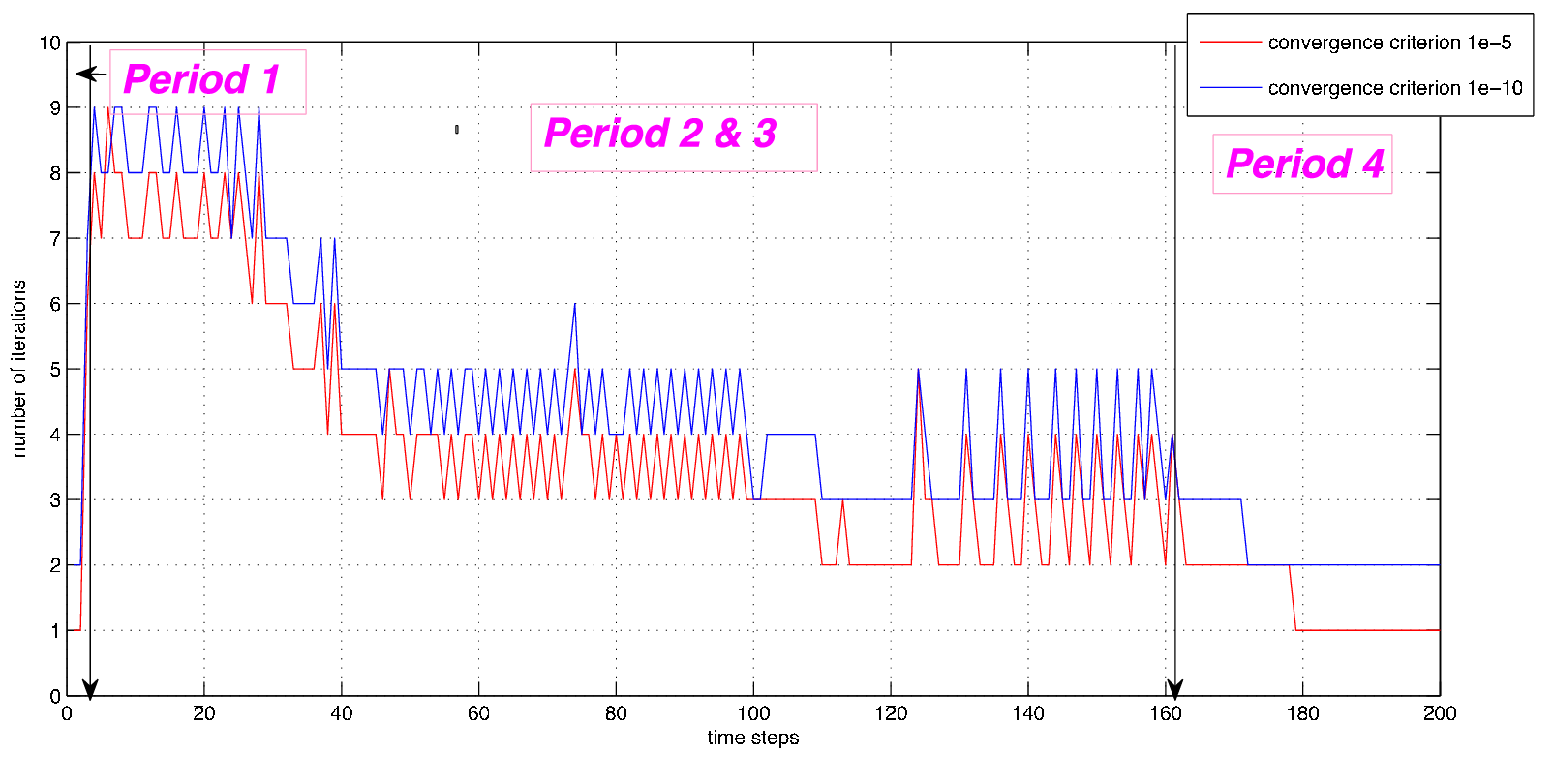}
\caption {\label{Quadratic}
Quadratic convergence of Newton-min: number of Newton-min iterations per time step
for two convergence criterions, 1.e-5 (red curve) and 1.e-10 (blue curve).}
\end{figure}
\end{center}

\medskip
The figure~\ref{Quadratic} shows the number of Newton-min iterations per time step for two convergence criterions, 
$\varepsilon_{1}=$ 1.e-5 (red curve) and $\varepsilon_{2}=$ 1.e-10 (blue curve). The points are connected with a straight line.
As mentioned at the end of section \ref{algorithm}, one can expect local quadratic convergence, at least for time steps which are sufficiently small.
In Figure~\ref{Quadratic}, we can observe this quadratic convergence.  Indeed one can verify in this figure  that, at each time step, the residue goes from 1.e-5 to 1.e-10 in one iteration. \
\newpage
\section{Conclusion}
We have studied a solution procedure for  a model describing
a system of two-phase (liquid-gas) flow in porous media with two components (hydrogen-water) where hydrogen 
can dissolve in the liquid phase. The problem is formulated as a nonlinear complementarity problem and is solved with the
Newton-min method. We considered an example of a Couplex-Gas benchmark and we showed 
 the ability of our solver to describe the appearance and disappearance of the gas phase during the
migration of hydrogen. We also discussed the quadratic convergence of the Newton-min method.
A theoretical justification for this quadratic convergence and other benchmark examples are under investigation.
\section*{Acknowledgments}
We thank J.\ Ch.\ Gilbert for discussions on complementarity problems and on the implementation of the Newton-min algorithm.










%

\end{document}